\documentclass{amsart}

\usepackage{amsmath,amssymb,amsfonts,graphicx}\usepackage[all]{xy}\usepackage{epstopdf}

\newtheorem{theorem}{Theorem}
\newcommand{\bt}{\begin{theorem}}
\newcommand{\et}{\end{theorem}}
\newtheorem{lemma}{Lemma}
\newcommand{\bl}{\begin{lemma}}
\newcommand{\el}{\end{lemma}}
\newtheorem{corollary}{Corollary}
\newcommand{\bc}{\begin{corollary}}
\newcommand{\ec}{\end{corollary}}
\newcommand{\bconj}{\begin{conjecture}}
\newcommand{\econj}{\end{conjecture}}
\newtheorem{problem}{Problem}
\newcommand{\bprob}{\begin{problem}}
\newcommand{\eprob}{\end{problem}}
\newcommand{\beq}{\begin{equation}}
\newcommand{\eeq}{\end{equation}}
\newcommand{\benum}{\begin{enumerate}}
\newcommand{\eenum}{\end{enumerate}}
\newcommand{\N}{\ensuremath{ \mathbf N }}
\newcommand{\Z}{\ensuremath{\mathbf Z}}
\newcommand{\Q}{\ensuremath{\mathbf Q}}
\newcommand{\R}{\ensuremath{\mathbf R}}

\newcommand{\mca}{\ensuremath{ \mathcal A}}

\newcommand{\mcm}{\ensuremath{ \mathcal M}}

\newcommand{\mcx}{\ensuremath{ \mathcal X}}

\newcommand{\mbc}{\ensuremath{ \mathbf c}}

\newcommand{\mbm}{\ensuremath{ \mathbf m}}

\newcommand{\mbo}{\ensuremath{ \mathbf 0}}
\newcommand{\mbu}{\ensuremath{ \mathbf u}}

\newcommand{\mbx}{\ensuremath{ \mathbf x}}

\newcommand{\bmat}{\left(\begin{matrix}}
\newcommand{\emat}{\end{matrix}\right)}
\newcommand{\bsmallmat}{\left(\begin{smallmatrix}}
\newcommand{\esmallmat}{\end{smallmatrix}\right)}

\DeclareMathOperator{\qqand}{\qquad\text{and}\qquad}

\date{\today}

\title[Sumset size races]{Sumset size races for measurable sets}

\author{Melvyn B. Nathanson}
\address{Department of Mathematics\\Lehman College (CUNY)\\Bronx, NY 10468}
\email{melvyn.nathanson@lehman.cuny.edu}

\date{\today}

\thanks{This work was supported in part by  PSC-CUNY Research Award Program grant 66197-00 54.}

\begin{document}

\begin{abstract}
Let $G$ be a locally compact abelian group with Haar measure $\mu$. 
For integers $n \geq 2$ and $H \geq 2$ and for 
any $n$-tuples $\mbu_1,\ldots, \mbu_H \in \N^n$, 
there exist measurable subsets 
$A_1,\ldots, A_n$ of $G$ such that the $n$-tuple 
$\left( \mu(hA_1),\ldots, \mu(hA_n) \right)$ has the same 
relative order as the $n$-tuple $\mbu_h$ for all $h = 1,\ldots, H$. 

For integers $m_{i,h}$ for $i =1,\ldots, n-1$ and $h  = 1,\ldots, H$, there are 
Lebesgue measurable sets $A_1,\ldots, A_n$ in \R\ such that $\mu(hA_{i+1}) - \mu(hA_i) = m_{i,h}$ for all $i$ and $h$. 
\end{abstract}

\date{\today}

\maketitle

\section{Sums of measurable sets} 

Let $\N = \{1,2,3,\ldots\}$ denote the set of positive integers 
and  $\N_0 = \{0,1,2,3,\ldots\}$ the set of nonnegative integers.  

Let $A$ be a nonempty subset of the integers or of any additive abelian group or 
semigroup.   
The $h$-fold sumset of $A$ is the set 
\[
hA = \{a_1 + \cdots + a_h: a_i \in A \text{ for all } i = 1,\ldots, h\}.  
\]  
Nathanson~\cite{nath26aa} introduced problems about ``sumset size races.'' 
Here is a simple example.  For every positive integer $m$, do there exist finite subsets $A$ and $B$ 
of the semigroup and sequences of positive integers $h_1 < h_2 < \cdots < h_m$ 
such that  the differences of the sumset sizes $|h_iA|$ and $|h_iB|$ 
oscillate for $i = 1,\ldots, m$, that is,
\begin{align*}
& |h_1A| < |h_1B|,  \qquad |h_2A| > |h_2B| \qquad |h_3A| < |h_3B|, \\
&  |h_4A| > |h_4B|,  \qquad |h_5A| < |h_5B| \qquad |h_6A| > |h_6B|, \ldots.
\end{align*}
Paul P\' eringuey and Anne de Roton solved this problem  for 
sumset sizes of finite sets of integers.   

\bt[P\' eringuey and de Roton~\cite{peri-roto25}]                   \label{races:theorem:roton}
For every integer $m \geq 3$, there exist finite sets $A$ and $B$ of integers 
and an increasing sequence of positive integers $h_1 < h_2 < \cdots < h_m$ 
such that $|A| = |B|$ and 
\[
|h_iA| < |h_iB| \qquad \text{if $i$ is odd}\\
\]
and 
\[
|h_iA| > |h_iB| \qquad \text{if $i$ is even.} 
\]
\et

The following more general problem was posed in~\cite{nath26aa}.  
The \emph{normalization} of an $n$-tuple $\mbu = (u_1,\ldots, u_n)$ 
of integers is the $n$-tuple $\tau(\mbu) = (\tau(1), \tau(2), \ldots, \tau(n))$ obtained 
by replacing the $i$th smallest entry in $\mbu$ with $i$ for all $i$.  
A \emph{$\tau$-tuple of length $n$} is an $n$-tuple of the form $\tau(\mbu)$ 
for some $n$-tuple $\mbu \in \Z^n$.
Here are examples of $6$-tuples $\mbu$ and their associated normalizations $\tau(\mbu)$:  
\vspace{0.3cm} 
\begin{center} 
\begin{tabular}{  l | c | c| c| c }  
$\mbu$         &    (1,2,3,4,5,6)   &  (-2,13,11,0,22,4)  &  (7,3,2,9,3,5)  & (9,7,8,9,7,8)  \\  \hline
$\tau(\mbu)$ &   (1,2,3,4,5,6)   & (1,5,4,2,6,3)        & (4,2,1,5,2,3) & (3,1,2,3,1,2) 
\end{tabular}
\end{center}
\vspace{0.3cm} 

Let $\mca = (A_1, A_2, \ldots, A_n)$ be an $n$-tuple of finite subsets 
of  an additive abelian semigroup.  
Associated to the $n$-tuple of $h$-fold sumsets $h\mca = (hA_1, hA_2, \ldots, hA_n)$ 
is the $n$-tuple of sumset sizes $|h\mca| =  ( |hA_1|, |hA_2|, \ldots, |hA_n|)$ 
and its normalization $\tau(|h\mca|)$.

\bprob[from~\cite{nath26aa}]        \label{races:problem:1}
Let $\mbu_1, \mbu_2, \ldots, \mbu_m$ be a finite sequence of $n$-tuples of integers.  
Does  there exist an $n$-tuple of  finite sets $\mca = (A_1, A_2, \ldots, A_n)$  
and a  sequence of positive integers 
$h_1 < h_2 < \cdots < h_m$ such that 
\[
\tau(|h_k \mca|) = \tau(|h_kA_1|, |h_kA_2|, \ldots, |h_kA_n|) = \tau(\mbu_k)
\]
for all $k = 1, \ldots, m$?
\eprob

Noah Kravitz solved  a strong form of this problem.  

\bt[Kravitz~\cite{krav25}]          \label{races:theorem:kravitz}
For integers $n \geq 2$ and $H \geq 2$ and for any $n$-tuples 
$\mbu_1,\ldots, \mbu_H,\mbu_{\infty}\in \N^n$, 
there exist finite subsets 
$A_1,\ldots, A_n \subseteq \Z$ such that 
\[
\left( \tau(|hA_1|,\ldots, |hA_n| \right) = \tau(\mbu_h)
\]
for all $h = 1,\ldots, H$ and 
\[
\tau(|hA_1|,\ldots, |hA_n|) = \tau(\mbu_{\infty})
\]
for all $h >  H$. 
\et

We prove a continuous analogue of this result. 

\bt         \label{races:theorem:Haar} 
Let $G$ be a locally compact  abelian group 
that is not a torsion group and let 
$\mu$ be a Haar measure on $G$.  
For integers $n \geq 2$ and $H \geq 2$ and for 
any $n$-tuples $\mbu_1,\ldots, \mbu_H \in \N^n$, 
there exist measurable subsets 
$A_1,\ldots, A_n$ of $G$ such that 
\[
\tau\left( \mu(hA_1),\ldots, \mu(hA_n) \right) 
= \tau(\mbu_h)
\]
for all $h = 1,\ldots, H$. 
\et

Theorem~\ref{races:theorem:Haar} 
solves Problem~\ref{races:problem:1} for Haar measure and, in particular, 
 for Lebesgue measure.

Jacob Fox, Noah Kravitz, and Shengtong  Zhang~\cite{FKZ25} obtained the following 
beautiful theorem, which is stronger than Theorem~\ref{races:theorem:kravitz} 
and proves much more than was asked in Problem 1.

\bt[Fox, Kravitz, and Zhang~\cite{FKZ25}]             \label{races:theorem:fox}
For integers $n \geq 2$ and $H \geq 2$, 
let $m_{i,h} \in \Z$ for all $i  = 1,\ldots, n-1 $ and $h  = 1,\ldots, H $.  
There exist finite sets of integers $A_1,\ldots, A_n$ such that 
\[
|hA_i| - |hA_{i+1}| = m_{i,h}
\]
 for all $i  = 1,\ldots, n-1 $ and $h  = 1,\ldots, H $. 
\et 

We prove a continuous analogue for subsets of \R.   

\bt         \label{races:theorem:nathanson-3} 
For integers $n \geq 2$ and $H \geq 2$, 
let $m_{i,h} \in \Z$ for all $i  = 1,\ldots, n-1 $ and $h  = 1,\ldots, H $.  
For all $\theta \in \R$, $\theta > 0$, there exist Lebesgue-measurable sets 
of real numbers $A_1,\ldots, A_n$ such that 
\[
\mu(hA_i) - \mu(hA_{i+1}) = \theta m_{i,h}
\]
for all $i  = 1,\ldots, n-1 $ and $h  = 1,\ldots, H $. 
\et

Note that, for all integers $j,k$ with $1 \leq j < k \leq n$, we have 
\[
\mu(hA_j) - \mu(hA_k)  = \sum_{i=j}^{k-1} \left( \mu(hA_i) - \mu(hA_{i+1}) \right)
\]
and so the $\binom{n}{2}$ differences $\mu(hA_j) - \mu(hA_k)$ are determined 
by the $n-1$ successive differences 
$\mu(hA_i) - \mu(hA_{i+1})$ for $i = 1, \ldots, n-1$.

\bprob
For integers $n\geq 2$ and $H \geq 2$, let $\theta_{i,h} \in \R$ 
for all $h = 1,2,\ldots, H$ and $i =1,2,\ldots, n-1$. 
Do there exist Lebesgue-measurable sets $A_1,\ldots, A_n$ of real numbers 
(or, in a locally compact additive abelian group, Haar-measurable sets 
$A_1,\ldots, A_n$) 
 such that 
\[
\mu(hA_i) - \mu(hA_{i+1}) = \theta_{i,h}
\]
for all $h = 1,2,\ldots, H$ and $i = 1,2,\ldots, n-1$? 
\eprob

It follows from  Theorem~\ref{races:theorem:nathanson-3} 
that the answer is ``yes'' if there   
exist integers $m_{i,j}$ and a real number $\theta$ 
such that $\theta_{i,h} = \theta m_{i,h}$ for all $i$ and $h$.  
What is the answer in the general case?

For related work on sumset sizes, see  Nathanson~\cite{nath25bb,nath25,nath25cc,nath26aa}, 
Rajagopal~\cite{raja25}, 
and Schinina~\cite{schi25}. 

\section{Proof of Theorem~\ref{races:theorem:Haar}}

The group $G$ is not a torsion group and so $G$ contains an element $x_0$ of 
infinite order.  Thus, for all $b,b' \in \Z$, we have $bx_0 = b'x_0$ if and only if $b = b'$.

By Kravitz's theorem, there are nonempty finite sets 
$B_1,\ldots, B_n$ of integers that satisfy 
\[
 \tau \left( |hB_1|,\ldots, |hB_n| \right)  = \tau(\mbu_h)
\]
for all $h = 1,\ldots, H$.  
The union of the sumsets $hB_i$ is the finite set  
\[
B =  \bigcup_{h=1}^H \bigcup_{i=1}^n hB_i. 
\]
Because a topological group is Hausdorff and $B$ is finite,  
$G$ contains an open neighborhood $V$ of the identity such that the open sets 
$bx_0+ V$ are pairwise disjoint for all $b \in B$.  
In a topological group, if $V$ is an open neighborhood of the identity, 
then for every positive integer $H$ there is an open neighborhood $U$ of the identity
such that $HU \subseteq V$.  It follows that  
\[
U \subseteq 2U \subseteq \cdots \subseteq HU \subseteq V.
\]
For $h = 1,\ldots, H$, the sumsets $hU$ are open and the 
open sets $bx_0+ hU$ are pairwise disjoint for $b \in B$.  
Nonempty open sets in  a topological group have positive Haar measure 
and so $\mu(bx_0+hU) = \mu(hU) >0$ for all $b \in B$. 

For all $i=1,\ldots, n$,  define the open set 
\[
A_i = \bigcup_{b_{i,j}\in B_i} \left( b_{i,j} x_0 + U \right). 
\]
For all $h = 1,\ldots, H$, we have the sumset 
\begin{align*}
hA_i 
& = h\left(  \bigcup_{b_{i,j}\in B_i} \left( b_{i,j} x_0 + U \right) \right) \\ 
& = \bigcup_{(b_{i,1},\ldots, b_{i,h}) \in B_i^h} \sum_{j=1}^h \left( b_{i,j} x_0 +U \right) \\
& = \bigcup_{(b_{i,1},\ldots, b_{i,h}) \in B_i^h} 
 \left(  \left(  \sum_{j=1}^h b_{i,j} \right)  x_0+ hU\right) \\
& =  \bigcup_{b\in hB_i} \left( bx_0+hU \right). 
\end{align*}
Because the sets $bx_0+ hU$ are pairwise disjoint, we have 
\begin{align*}
\mu\left( hA_i \right) 
& = \mu\left(  \bigcup_{b\in hB_i} \left( bx_0+ hU \right) \right) \\ 
& = \sum_{b\in hB_i}  \mu\left(  bx_0+ hU \right) \\
& = |hB_i| \mu(hU)
\end{align*}
for all $i = 1,\ldots, n$.  Because $\mu(hU)>0$, we obtain 
\begin{align*}
 \tau\left( \mu(hA_1),\ldots, \mu(hA_n) \right) 
& = \tau \left(  |hB_1| \mu(hU),\ldots, |hB_n| \mu(hU) \right) \\
& = \tau  \left( |hB_1|,\ldots, |hB_n| \right) \\
& = \tau(\mbu_h). 
\end{align*} 
This completes the proof.

\section{Proof of Theorem~\ref{races:theorem:nathanson-3}}     \label{races:section:proof} 
This proof closely follows the Fox-Kravitz-Zhang 
proof of Theorem~\ref{races:theorem:fox} for sets of integers. 
  
Let $A$ be a set of real numbers.  
The dilation of the set $A$ by the real number $\lambda$ is the set 
$\lambda \ast A = \{\lambda a:a \in A\}$. 
If $A$ is Lebesgue-measurable, then $\mu(\lambda\ast A) = |\lambda|\mu(A)$.

\bl               \label{races:lemma:1}
Let $0 < \varepsilon \leq 1/3$ and 
\[
X = [0,1-\varepsilon] \cup [2,3-\varepsilon]. 
\]
For all integers $h \geq 3$,
\[
hX = [0,h(3-\varepsilon)].
\]
If $Y$ is a nonempty subset of  $[1+\varepsilon, 2-2\varepsilon]$, then 
\[
2(X \cup Y) = [0,2(3-\varepsilon)].
\]
\el

\begin{proof}
We have 
\begin{align*}
hX  & = h ( [0,1-\varepsilon] \cup [2,3-\varepsilon] ) \\
& = \bigcup_{k=0}^h \left( k  [0,1-\varepsilon] + (h-k) [2,3-\varepsilon]  \right)        \\ 
& = \bigcup_{k=0}^h \left(    [0,k(1-\varepsilon)] + [2 (h-k),  (h-k)(3-\varepsilon)]  \right)  \\ 
& =  \bigcup_{r=0}^h \left[ 2(h-k),  k(1-\varepsilon) +  (h-k) (3-\varepsilon) \right]     \\ 
& =  \bigcup_{r=0}^h \left[ 2h - 2k,  h(3-\varepsilon)  -2 k\right].  
\end{align*}
These $h+1$ intervals ``move to the left'' as $k$ increases and overlap if
\[
 2h - 2(k-1) \leq h(3-\varepsilon)  -2 k 
\]
or, equivalently, if 
\[
2 \leq h(1-\varepsilon).  
\]
If  $\varepsilon \leq 1/3$, then this inequality holds for all $h \geq 3$ and 
\[
hX = [0,h(3 - \varepsilon)].  
\]

For $h = 2$, we have 
\[
2X = [0,2-2\varepsilon] \cup [2, 4-2\varepsilon] \cup [4, 6 - 2\varepsilon].
\]
If $y \in Y \cap [1+\varepsilon, 2-2\varepsilon]$, then 
\[
y \leq 2 - 2\varepsilon < 2 = (1-\varepsilon) + (1 + \varepsilon) \leq    1-\varepsilon  + y 
\]
and so
\[
[2-2\varepsilon, 2] \subseteq [ y  , 1-\varepsilon + y] = [ 0, 1-\varepsilon] + \{y\}  
\subseteq X+ \{y\} \subseteq X+ Y.
\]
Similarly,
\[
2+y \leq 4-2\varepsilon < 4 = (3-\varepsilon) + (1+\varepsilon) \leq 3-\varepsilon + y 
\]
and so 
\[
[ 4- 2\varepsilon, 4] \subseteq [2+y , 3-\varepsilon + y] =  [ 2, 3-\varepsilon] + \{y\}   
\subseteq X + \{y\}  \subseteq X+ Y.  
\]
It follows that  $2(X +Y) = [0,2(3-\varepsilon)]$. 
This completes the proof. 
\end{proof}

\bl                       \label{races:lemma:2}
For $H \geq 2$, let $\varepsilon$, $\delta$, and $c$ 
be real numbers such that    
\beq              \label{races:inequalities} 
0 < \varepsilon <\frac{1}{3}, \qquad 0 < \delta < \frac{\varepsilon}{H-1}, 
\qquad  (H-1) \delta + 3 < c.  
\eeq
Let 
\[
X = [0,1-\varepsilon] \cup [2,3-\varepsilon]  
\]
and, for all $j = 1,\ldots, n$,  let   
\[
Y_j \subseteq [1+\varepsilon, 2-2\varepsilon] 
\]
be a nonempty measurable set.  Let 
\[
A_j = [0,\delta] \cup (c+ (X \cup   Y_j) ).
\]
For all integers $h   \in [1,H]$ and $j,k \in [1,n]$,  
\beq          \label{races:2X}
\mu(hA_j) - \mu(hA_k)  
=  \mu\left( [0, (h-1) \delta] + Y_j \right) -  \mu\left(  [0, (h-1) \delta] + Y_k  \right). 
\eeq
\el

\begin{proof} 
For all $h = 1,2,\ldots, H$ and $j = 1,\ldots, n$, we have 
\begin{align*}
hA_j & = h \left( [0,\delta] \cup (c+ (X \cup   Y_j) )\right) \\ 
& = \bigcup_{r=0}^h  \left( (h-r)[0,\delta]  + r(c+(X\cup Y_j)) \right)    \\ 
& =   [0, h\delta] \cup \left(  [0, (h-1) \delta]  +c +  (X\cup Y_j) \right)    \\
& \qquad  \cup  \bigcup_{r=2}^h \left( [0, (h-r)\delta] + rc + r(X\cup Y_j) \right).
\end{align*}
Applying  Lemma~\ref{races:lemma:1}, we obtain 
\[
 r(X\cup Y_j) =  [0,r(3-\varepsilon)]  
\]
for all $r \geq 3$ and so 
\begin{align}
hA_j 
& =   [0, h\delta] \cup \left(  [0, (h-1) \delta]  +c +  (X\cup Y_j) \right)  \label{races:HAj}   \\
& \qquad  \cup  \bigcup_{r=2}^h \left( [0, (h-r)\delta] + rc + [0,r(3-\varepsilon)] \right).  \nonumber
\end{align}

Inequalities~\eqref{races:inequalities} imply the following inequalities: 
\begin{align*}
\max\left(  [0, h\delta]  \right) & = h\delta \leq H\delta < (H-1)\delta +3 \\
& < c  = \min \left( [0, (h-1) \delta]  +c +  (X\cup Y_j) \right) 
\end{align*} 
and
\begin{align*}
\max ( [0, (h- & 1) \delta]  +  c + (X\cup Y_j)   ) = (h-1)\delta + c + 3-\varepsilon \\ 
& <  2c = \min\left(  \bigcup_{h=2}^h \left( [0, (h-h)\delta] + hc + [0,h(3-\varepsilon)] \right)  \right). 
\end{align*}
It follows that the three sets in the union~\eqref{races:HAj}  are pairwise disjoint. 
For $j,k \in [1,n]$, translation invariance of Lebesgue measure implies  
\begin{align*} 
\mu & \left( hA_j \right)  - \mu\left(hA_k \right) \\ 
&\qquad = \mu \left(  [0, (h-1) \delta]  +c +  (X\cup Y_j) \right)  -  \mu\left( [0, (h-1) \delta]  +c +  (X\cup Y_k)  \right)  \\ 
& \qquad =  \mu\left(  [0, (h-1) \delta]  +  (X\cup Y_j) \right)  -  \mu\left( [0, (h-1) \delta]  +  (X\cup Y_k)  \right)  .
\end{align*} 

Because 
\[
X = [0,1-\varepsilon] \cup [2,3-\varepsilon]  
\]
we obtain  
 \begin{align}
 [0, (h-1) \delta] +  (X\cup Y_j)  
= & \left( [0, (h-1) \delta]   +  [0,1-\varepsilon] \right)   \nonumber   \\
&\cup \left( [0, (h-1) \delta]   +  Y_j  \right)  \label{races:XUY}   \\
&  \cup \left( [0, (h-1) \delta]   + [2,3-\varepsilon]   \right).   \nonumber
\end{align}
Inequalities~\eqref{races:inequalities} imply the following inequalities: 
 \begin{align*}
\max & \left(  [0, (h-1) \delta]   +  [0,1-\varepsilon]  \right) 
=  (h-1) \delta    + 1-\varepsilon \\
&  \leq (H-1)\delta +1 - \varepsilon  <  1 <  1+\varepsilon \\
& \leq \min \left([0, (h-1) \delta]   +  Y_j    \right) 
\end{align*}
and
 \begin{align*}
\max &\left( [0, (h-1) \delta]   +  Y_j   \right) 
  \leq  (h-1) \delta   + 2 - 2\varepsilon  \\
&  \leq  (H-1) \delta   + 2 - 2\varepsilon  < 2 - \varepsilon  <  2 \\
&  = \min \left(  [0, (h-1) \delta]   + [2,3-\varepsilon]  \right).  
\end{align*}
It follows that the three sets in the union~\eqref{races:XUY} are pairwise disjoint, and so 
\begin{align*}
 \mu  \left( hA_j\right) -  \mu\left(  hA_k \right) & =  \mu\left(  [0, (h-1) \delta]  +  (X\cup Y_j) \right) - \mu\left( [0, (h-1) \delta]  +  (X\cup Y_k)  \right)  \\
& = \mu\left( [0, (h-1)\delta] + Y_j  \right)-\mu\left( [0, (h-1) \delta] + Y_k \right).
\end{align*} 
This completes the proof. 
\end{proof}

\bl                       \label{races:lemma:3}
For $H \geq 2$, let $\ell_1,\ell_2,\hdots, \ell_H$ 
be nonnegative integers, not all 0, and  let $ L_H = \sum_{r=1}^H \ell_r$. 
Let 
\beq         \label{races:inequalities-2}
0 < \varepsilon < \frac{1}{3} 
\qqand 
0 < \delta \leq \frac{1-3\varepsilon}{2HL_H}. 
\eeq 
The closed interval $[1+\varepsilon, 2 - 2\varepsilon]$ contains a set $Z$ that is the 
union of $L_H$ pairwise disjoint open intervals  
with the property that the nonempty set 
\[
Y = [1+\varepsilon, 2 - 2\varepsilon] \setminus Z 
\]
satisfies 
\[
\mu \left( [0,(h-1)\delta] + Y \right)  
=  (h-1)\delta  + 1-3\varepsilon - \delta \sum_{r=h}^H (r-h+1) \ell_r. 
\]
for all $h = 1, \ldots, H$. 
\el

\begin{proof}
The interval $[1-\varepsilon, 2-2\varepsilon]$ has length $1-3\varepsilon$. 
In this interval we construct $L_H$ pairwise disjoint open subintervals,  
each of length $2H\delta$, where 
\[
0 < \delta \leq \frac{1-3\varepsilon}{2HL_H}.
\]
For $j = 0, 1,2,\ldots,L_H$, let 
\[
u_j =  1 + \varepsilon + 2Hj\delta.        
\]
Then $u_j - u_{j-1} = 2H\delta$ for all $j=1,\ldots, L_H$ and 
\begin{align*}
1+\varepsilon = u_0 
&< u_1 < \cdots < u_{L_H-1} \\
& < u_{L_H}
=  1 + \varepsilon + 2HL_H \delta \\
& \leq 2-2\varepsilon.  
\end{align*}
Let $L_0 = 0$ and 
\[
L_r = \sum_{i=1}^r \ell_i 
\] 
for $r = 1,\ldots, H$.  
We construct $L_H$ open  intervals $I_j$ as follows.  
Every integer $j \in [1,L_H]$  can be written uniquely in the form 
\[
j = L_{r-1}+i
\]
for some $r \in \{1,\ldots,H\}$ and $i \in \{ 1, \ldots, \ell_r \}$. Let
\[
I_j = (u_j, u_j+r\delta) 
\]
for all  $j = 1,2,\ldots,L_H$.  
For all $r = 1,\ldots, H$ there are $\ell_r$ intervals $I_j$ of length $r\delta$. 
The intervals $I_j$ are pairwise disjoint because 
\[
u_j+r\delta \leq u_j + H\delta < u_j+2H\delta = u_{j+1} 
\]
for all $j = 1,\ldots, L_H-1$.  This inequality also implies 
\[
u_j + [0,(h-1)\delta] =  [u_j, u_j + (h-1)\delta] \subseteq [u_j, u_{j+1}]  
\]
for all $h = 1,\ldots, H$.  

Let 
\[
Z = \bigcup_{j=1}^{L_H}  I_j \qqand Y = [1+\varepsilon,2-2\varepsilon]\setminus Z. 
\]
The set $Y$ is nonempty with  $u_j \in Y$ for all $j = 1, \ldots, L_H$. 
We have 
\[
\mu(Z) = \sum_{j=1}^{L_H} \mu(I_j) = \delta \sum_{r=1}^H r \ell_r 
\]
and 
\begin{align*} 
\mu(Y) & = 1-3\varepsilon - \mu(Z)  \\
& = 1-3\varepsilon  -  \delta \sum_{r=1}^H r \ell_r. 
\end{align*}

Let $h \in \{1,\ldots, H \}$.  We shall  compute 
\[
 \mu\left(  [0, (h-1) \delta] + Y \right).    
\] 
There are three steps.  
First, we observe that $ 2-2\varepsilon \in Y$ implies 
\begin{align*} 
 [2-2 \varepsilon ,2-2 \varepsilon + (h-1)\delta]  
 & = 2-2\varepsilon   + [0,(h-1)\delta] \\
& \subseteq [0,(h-1)\delta] + Y. 
\end{align*}

Second, if $1\leq j \leq L_{h-1}$, then $L_{r-1} < j \leq L_r$ 
for some $r \in [1,h-1]$ and 
\begin{align*} 
I_j & = (u_j, u_j+r\delta) \subseteq (u_j, u_j+ (h-1)\delta) \\
& = u_j + (0,(h-1)\delta) \\
& \subseteq [0,(h-1)\delta] + Y.  
\end{align*}
Therefore, 
\[
 I_j \setminus  ( [0,(h-1)\delta] +Y) = \emptyset 
\]
and 
\[
\mu\left(  I_j \setminus  ( [0,(h-1)\delta] +Y)  \right) = 0.
\]

Third, if $L_{h-1} < j \leq L_H$, then $L_{r-1} < j \leq L_r$ 
for some $r \in [h, H]$ and  
\[ 
I_j   = (u_j, u_j+r\delta)   \supseteq (u_j, u_j+ (h-1)\delta].
\]
It follows that 
\begin{align*}
I_j  \setminus (u_j, u_j+ (h-1)\delta]  = (u_j+ (h-1)\delta, u_j + r\delta ) 
\end{align*}
and 
\[
\mu\left(I_j \setminus ((0,(h-1)\delta) +Y ) \right) = (r-h+1)\delta. 
\]

Combining the three steps, we obtain 
\begin{align*} 
[1+\varepsilon, (h-1)\delta + & 2-2\varepsilon] \setminus  [0,(h-1)\delta] + Y \\
 & = \left(
 \bigcup_{r=h}^{H} \bigcup_{j=L_{r-1}+1}^{L_r} (u_j + (h-1)\delta, u_j+r\delta) 
 \right) 
\end{align*} 
The intervals $ (u_j + (h-1)\delta, u_j+r\delta) $ are pairwise disjoint and so 
\begin{align*} 
\mu & \left(
 \bigcup_{r=h}^{H} \  \bigcup_{j=L_{r-1}+1}^{L_r} (u_j + (h-1)\delta, u_j+r\delta) 
 \right) \\
 &  \qquad=   \sum_{r=h}^{H} \ \sum_{j=L_{r-1}+1}^{L_r} 
 \mu (u_j + (h-1)\delta, u_j+r\delta) \\
 & \qquad =  \delta \sum_{r=h}^{H} (r-h+1) \ell_r .
 \end{align*} 
Therefore, 
\[
\mu \left( [0,(h-1)\delta] + Y  \right) = (h-1)\delta + 1-3\varepsilon 
-  \delta \sum_{r=h}^{H} (r-h+1) \ell_r . 
\]
This completes the proof. 
\end{proof}

\bl                       \label{races:lemma:4}
For integers $H \geq 2$ and $n \geq 2$ and 
for $i  = 1,\ldots, n-1 $ and $j =1,2,\ldots, H $, 
let  $\ell_{i,j}$ be nonnegative integers such that 
\[
L_{i,H} =   \sum_{j=1}^H \ell_{i,j} > 0
\]
for all $i$.  Let $\varepsilon$,  $\delta$, and $c$ be real numbers such that 
\[
0  < \varepsilon < \frac{1}{3} 
\]
and 
\[
0 < \delta< \min\left\{ \frac{\varepsilon}{H-1},\frac{1-3\varepsilon}{2HL_{i,H}} \right\}.  
\]  
and 
\[
(H-1)\delta + 3 < c.
\]
The  interval $[1+\varepsilon, 2 - 2\varepsilon]$ contains nonempty 
measurable sets $Y_1,\ldots, Y_n$ such that the sets $A_1,\ldots, A_N$ defined by 
\[
A_i = [0,\delta] \cup (c+(X+Y_i)) 
\]
satisfy 
\[
\mu(hA_i) - \mu(hA_{i+1}) 
=   \delta \sum_{r=h}^H (r-h+1) (\ell_{i+1,r} -  \ell_{i,r} ) 
\]
for all $h  = 1, \ldots, H$. 
\el

\begin{proof}
Fom Lemmas~\ref{races:lemma:2} and~\ref{races:lemma:3} we obtain 
\begin{align*}
\mu(hA_i) - \mu(hA_{i+1}) 
& = \mu \left(  [0,(h-1)\delta] + Y_i \right)  - \mu \left(  [0,(h-1)\delta] + Y_{i+1} \right) \\ 
& = \left( (h-1)\delta  + 1-3\varepsilon - \delta \sum_{r=h}^H (r-h+1) \ell_{i,r} \right) \\
& \qquad -  \left( (h-1)\delta  + 1-3\varepsilon - \delta \sum_{r=h}^H (r-h+1) \ell_{i+1,r}  \right) \\
& =   \delta \sum_{r=h}^H (r-h+1) (\ell_{i+1,r} -  \ell_{i,r} ) 
\end{align*} 
for all $h  = 1, \ldots, H$. 
\end{proof}

The last piece of the proof of Theorem~\ref{races:theorem:nathanson-3} 
is the following result about linear Diophantine equations.

\bl                       \label{races:lemma:5}
For integers $H \geq 2$ and $n \geq 2$, let  $m_{i,j}  \in \Z$
 for all $i  = 1,\ldots, n-1 $ and $j =1,2,\ldots, H$. 
There exist $\ell_{i,j} \in \N_0$ 
 for all $i =1,\ldots, n $ and $j =1,2,\ldots, H$  
such that   
\beq            \label{races:upperTriangular}
\sum_{r=h}^H (r-h+1) (\mathbf{\ell}_{i+1,r} - \mathbf{\ell}_{i,r}) = m_{i,h} 
\eeq 
 for all $i  = 1,\ldots, n-1 $ and $ h =1,2,\ldots, H$.   
\el

\begin{proof}
For all $i  = 1,\ldots, n-1$, 
consider the following system of $H$ linear equations 
in $H$ variables $x_{i,1},\ldots, x_{i,H}$:
\[
\sum_{r=h}^H (r-h+1) x_{i,r} = m_{i,h} 
\]
for $h =1,\ldots, H$.
 The matrix for this system of equations is upper triangular with 
integer coefficients and 1's on the main diagonal: 
\[
M = \bmat 1 & 2 & 3 &4 &  \cdots & H-1 & H \\
0 & 1 & 2 & 3 & \cdots & H-2 & H-1 \\ 
0 & 0 & 1 & 2   & \cdots & H-3 & H-2 \\
\vdots & &&&&& \vdots \\
0 & 0 & 0 & 0 & \cdots & 1 & 2 \\
0 & 0 & 0 & 0 & \cdots & 0 & 1
\emat 
\]
and so the system has a unique solution in integers $x_{i,1},\ldots, x_{i,H}$.  
The explicit solution is 
\begin{align*}
x_{i,H} & = m_{i,H}  \\
x_{i,H-1} & = m_{i,H-1} - 2m_{i,H} \\ 
x_{i,r} & =  m_{i,r} -2m_{i,r+1}+m_{i,r+2} \qquad\text{for all $ r \in [1, H-2]$.}
\end{align*}

 Consider next, for $r=1,\ldots, H$,  the linear system of $n-1$ equations 
 in $n$ variables $\ell_{i,r}$: 
\begin{align*}
\ell_{2 ,r}   - \ell_{ 1,r}  & = x_{ 1,r} \\
\ell_{ 3,r}   - \ell_{ 2,r}   & = x_{2 ,r}  \\
\ell_{ 4,r}   - \ell_{ 3,r}  & = x_{3 ,r} \\
\vdots & \\
\ell_{ i+1,r}  - \ell_{i ,r}   & = x_{ i,r} \\
\vdots & \\
\ell_{n ,r}  - \ell_{n-1 ,r}   & = x_{n-1 ,r}
\end{align*}
This system has a unique solution for every choice of $\ell_{1,r}$: 
\[
 \ell_{i,r} = \ell_{ 1,r}  + \sum_{j=1}^{i-1} x_{j,r} \qquad\text{for $i=2,\ldots, n$.}
 \]
Choosing 
\[
\ell_{ 1,r}  \geq \max\left( 0, -\sum_{j=1}^{i-1} x_{j,r} : i = 1,\ldots, n-1 \right)
\]
gives a nonnegative solution. 
If the $x_{i,r}$ are integers, then choosing 
\[
\ell_{ 1,r} = \max\left( 0, -\sum_{j=1}^{i-1} x_{j,r} : i  = 1,\ldots, n-1 \right)
\]
gives the smallest nonnegative integral solution.  
This completes the proof. 
\end{proof}

We can now complete the proof of Theorem~\ref{races:theorem:nathanson-3}.
Given integers $m_{i,h}$ for all $i = 1,\ldots, n-1$ 
and $h = 1,\ldots, H$, we apply Lemmas~\ref{races:lemma:4} 
and~\ref{races:lemma:5} and obtain  integers $\ell_{i,r}$ 
and measurable sets $A_1,\ldots, A_n$ such that 
\begin{align*}
\mu(hA_i) - \mu(hA_{i+1})  
& =  \delta \sum_{r=h}^H (r-h+1) (\ell_{i+1,r} -  \ell_{i,r} ) \\
& = \delta m_{i,h}.
\end{align*} 
Dilating the sets $A_1,\ldots, A_n$ by $\theta/\delta$ 
completes the proof.


\begin{thebibliography}{99} 


\bibitem{FKZ25} 
J. Fox, N. Kravitz and S. Zhang, 
Finer control on relative sizes of iterated sumsets,  
arXiv: 2506.05691.


\bibitem{krav25}
N. Kravitz, Relative sizes of iterated sumsets, 
J. Number Theory 272 (2025), 113--128, 



\bibitem{nath25bb}
M.  B. Nathanson, Problems in additive number theory, VI:  
Sizes of sumsets of finite sets, 
Acta Math. Hungar. 176 (2025), 498--521. 
doi.org/10.1007/s10474-025-01559-7 


\bibitem{nath25cc}
M.  B. Nathanson, 
Explicit sumset sizes in additive number theory, 
Canad. Math.  Bull. (2025), 1--12. 
doi.org/10.4153/S0008439525101549

\bibitem{nath25}
M.  B. Nathanson, 
Compression and complexity for sumset sizes in additive number theory,
J. Number Theory 281 (2026), 321--343.   
doi.org/10.1016/j.jnt.2025.09.025 



\bibitem{nath26aa}
M.  B. Nathanson, Inverse problems for sumset sizes of finite sets of integers, 
Fibonacci Quar. (2026), to appear. 
arXiv:2411.02365.

\bibitem{peri-roto25}
P. P\' eringuey and A. de Roton, 
A note on iterated sumsets races,
arXiv:2505.11233.  


\bibitem{raja25}
I. Rajagopal, 
Possible sizes of sumsets, 
arXiv: 2510.23022.  


\bibitem{schi25}
V. Schinina, On the sumset  of sets of size $k$, 
Integers, to appear. 
arXiv:2505.07679.  


\end{thebibliography}
\end{document}